\documentclass[letterpaper]{article}

\usepackage{fullpage,authblk}
\usepackage{amssymb,amsfonts,amsmath,enumerate}
\usepackage[colorlinks=false]{hyperref}
\usepackage{graphicx,subfigure,bm,latexsym}
\usepackage{multirow}
\usepackage{appendix}
\usepackage{float}
\usepackage{lineno}
\usepackage{natbib}
\setcitestyle{authoryear,open={(},close={)}}

\begin{document}



\title{ERFit: Entropic Regression Fit Matlab Package, for Data-Driven System Identification of Underlying Dynamic Equations}

\author[1,2]{Abd AlRahman AlMomani}
\author[1,2]{Erik Bollt}

\affil[1]{Department of Electrical and Computer Engineering, Clarkson University, Potsdam, NY 13699, USA}
\affil[2]{Clarkson Center for Complex Systems Science ($C^3S^2$), Potsdam, NY 13699, USA}

\date{}
\maketitle

\begin{abstract}
    Data-driven sparse system identification becomes the general framework for a wide range of problems in science and engineering. It is a problem of growing importance in applied machine learning and artificial intelligence algorithms. In this work, we developed the Entropic Regression Software Package (ERFit), a MATLAB package for sparse system identification using the entropic regression method. The code requires minimal supervision, with a wide range of options that make it adapt easily to different problems in science and engineering. The ERFit is available at \url{https://github.com/almomaa/ERFit-Package}
\end{abstract}








\subsection*{Keywords}

Data-Driven, System Identification, Machine Learning, Mathematical Modeling, Nonlinear Dynamics, Sparse Regression, Information Theory, Causality Inference, Entropic Regression.

\section{Introduction}

In this paper, we introduce an efficient Matlab implementation for the Entropic Regression method we introduced in \cite{almomani2020entropic,almomani2019prediction}. Entropic Regression (ER), is a data-driven discovery method for the underlying dynamics using sparse system identification. ER uses the conditional mutual information as an information-theoretic criterion and iteratively select relevant basis functions in a greedy search optimization scheme in terms of information criterion objective. Consider the problem in the matrix form:

\begin{equation} \label{eq:mainMatrixForm}
    \dot{X} = F(X) =  \Phi(X) \mathbf{\beta} 
\end{equation}
where $X\in\mathbb{R}^{N\times d}$ is the measured state variables of the $d$-dimensional system with $N$ observations, $\dot{X}$ is the vector field estimated from $X$, $\Phi:\mathbb{R}^{N\times d} \mapsto \mathbb{R}^{N\times K}$ is a function that maps the state variables $X$, to expanded set of candidate functions (this general form in Eq.~\ref{eq:mainMatrixForm}, cores a wide range ordinary differential equations, writing the vector field as a linear combination of possibly nonlinear basis functions) \cite{carleman1932}, and $\mathbf{\beta}\in\mathbb{R}^{K\times d}$ is the parameters matrix. For the sake of clarity we will write $\Phi(X)$ as $\Phi$ in the following discussion.

The basis functions do not need to be mutually orthogonal, and this was the main theme in previous approaches on nonlinear SID, with different methods differ mainly on how a model's fit is quantified \cite{crutchfield1987equations, hamilton2015predicting}. The different approaches include using standard squared error measures~\cite{BolltYao2007,Chen1989}, sparsity-promoting methods~\cite{Kalouptsidis2011,Brunton2015,Wang2011,Wang2016,kaiser2018sparse, brunton2016sparse, TranWard} as well as using entropy-based cost functions~\cite{guo2008extended}. 
Among those, sparsity-promoting methods have proven particularly useful because they tend to avoid the issue of overfitting, thus allowing a large number of basis functions to be included to capture possibly rich dynamical behavior~\cite{Kalouptsidis2011,Brunton2015,Wang2011}.

The inverse problem is then: Given $\Phi$ and $\dot{X}$, find $\beta$. The Best Linear Unbiased Estimator (BLUE) \cite{henderson1975best}, of the parameters matrix is known to be the least squares solution given by:
\begin{eqnarray}
\mathcal{L}(\dot{X},\Phi) & = & (\Phi^T\Phi)^{-1}\Phi^{T}\dot{X} \notag \\ 
                         & = & \Phi^{\dagger}\dot{X}
\end{eqnarray}
where $\Phi^{\dagger}$ is the pseudoinverse of the matrix $\Phi$. The reconstructed vector field using the least squares solution is given by:
\begin{eqnarray}\label{eq:LSprojection}
\mathcal{V}(\dot{X},\Phi) & = &  \Phi\Phi^{\dagger}\dot{X} \notag \\
                         & = &  \Phi\mathcal{L}(\dot{X},\Phi)
\end{eqnarray}

Sparse Regression problem now finding the minimal set of index $s\subset \{1,2,...,K\}$, such that $\mathcal{V}(\dot{X},\Phi_{s})$ as close as possible for $\dot{X}$ according to adopted quality measure (or objective function, cost function, loss function), where $\Phi_{s} \subset \Phi$ is the matrix with only the columns with index $i\in s$.

The ER method is a greedy search optimization method that contains two stages: Forward ER and Backward ER; in both stages, selection and elimination of basis functions are based on an entropy criterion (conditional mutual information). 

\subsection*{Forward Selection}
In the forward stage, our objective is to select the subset $s\subset \mathcal{S}=\{1,2,...,K\}$, that represent a strong candidate functions. Starting from empty set $s_0 =\emptyset$, the forward selection stage can be written as:
\begin{eqnarray}\label{eq:forward}
u_k & = & arg\displaystyle\max_{i\in\mathcal{S},i\notin s_{k-1}} I(\dot{X}; \mathcal{V}(\dot{X},\Phi_i)|\mathcal{V}(\dot{X},\Phi_{s_{k-1}})), \notag \\
s_k & = & s_{k-1} + u_k
\end{eqnarray}
where $k=1,...$, is the iteration index, $u_k$ is the set of index with the maximum objective function value. Note that $s_0=\emptyset \implies \mathcal{V}(\dot{X},\Phi_{s_0}) =\emptyset$ which reduces the conditional mutual information $I(\cdot; \cdot|\cdot)$ to the mutual information $I(\cdot;\cdot)$. The forward stage have a reward function, where at each iteration $k$, given the information ($\mathcal{V}(\dot{X},\Phi_{s_{k-1}})$) we already have from the set $s_{k-1}$, we are looking for the function that maximally add extra information to the model. The process terminates when the termination condition $HLT1$ met, which we will describe in the following sections.

\subsection*{Backward Elimination}
After the termination of the forward selection, we have the set $s$ that has the indices of the strong candidate functions. Eventually, $s$ may have a few non-relevant functions indices that are selected due to a high degree of uncertainty and the rounding error at the end of forward ER. Since we have reduced set of functions indices ($card(s)<<K$), it would be inexpensive to perform a validation operation to ensure the accuracy of the model, and the backward ER represent this operation. The backward stage is an elimination stage, where the functions indexed by $s$ re-examined for their information-theoretic relevance and these that are redundant will be removed. In particular, we label the set $s$ as initial set $s_0 = s$ for the backward stage, and we perform the following computations and updates,

\begin{eqnarray}\label{eq:backward}
u_k & = & arg\displaystyle\min_{i\in s_{k-1}} I(\dot{X}; \mathcal{V}(\dot{X},\Phi_i)|\mathcal{V}(\dot{X},\Phi_{\{s_{k-1}-i\}})), \notag \\
s_k & = & s_{k-1} - u_k.
\end{eqnarray}
The backward stage has a loss function, where at each iteration $k$, we examine information that will potentially be lost if we remove the index $i$ from the set $s_{k-1}$. We continue the elimination process until the termination condition $HLT2$ met.
 
The result of the backward ER is a set of indices $s$. We emphasize here that the forward ER stage can greatly reduce the computational complexity of the backward stage, by limiting the elimination search space to a few candidate functions. However, the backward elimination has a lower rate of error than the forward selection, and in case we having a low-dimension system, or we have efficient computations resources, we can skip the forward stage, and apply the backward stage directly with initial set $s = \{1,2,\dots,K\}$, and we provide this option in our software package.

The corresponding parameters $\mathbf{\beta} \in\mathbb{R}^{K}$, can be found by updating the vector of zeros $\mathbf{\beta}=\mathbf{0}_K$ such that:
\begin{equation}
    \mathbf{\beta}_s = \mathcal{L}(\dot{X},\Phi_s)
\end{equation}
where $\mathbf{\beta}_s$ is the entries of $\mathbf{\beta}$ indexed by the elements of $s$, which gives that $\Vert \mathbf{\beta} \Vert_0 = card(s)$. Note that the ER focus on finding the optimal set of basis functions in terms of the conditional mutual information with the vector field, and after finding this set, we find the value of the parameters by the ordinary least squares, without any attempts to apply any advanced techniques for optimizing the parameter's magnitude. In our package, we provide the function (\texttt{getSystemHandle.m}), which translate the estimated parameters to a system of ODEs, and return a function handle of the ODE function.

\subsection*{Termination Conditions}
In theory, the mutual information $I(x;y|z)$ is always non-negative and equals zero if and only if $x$ and $y$ are statistically independent given $z$. However, in practice, due to finite sampling and estimation inaccuracies, the estimated mutual information does not always equal to zero even when $x$ and $y$ are independent and can be negative. Thus, one needs a way to decide whether $x$ and $y$ should be deemed independent given the estimated value of $I(x;y|z)$. In \cite{sun2014identifying}, the authors introduced a standard shuffle test, with a ``confidence'' parameter $\alpha\in[0,1]$ for tolerance estimation. The shuffle test requires random shuffling of one of the variables $n_s$ times, to build a test statistic. In particular, for the $i$-th random shuffle, a random permutation $\pi^{(i)}:[T]\rightarrow[T]$ is generated to shuffle one of the variables, say $y$, which produces a new variable $(\tilde{y}^{(i)})$ where $\tilde{y}^{(i)}=y_{\pi^{(i)}}$; $x$ and $z$ are kept the same. Then, we estimate the mutual information $I(x;\tilde{y}^{(i)}|z)$ using the (partially) permuted variable $(x,\tilde{y}^{(i)},z)$, for each $i=1,\dots,n_s$. For given $\alpha$, we then compute a threshold value $I_\alpha(x;y|z)$ as the $\alpha$-percentile from the values of $I(x;\tilde{y}^{(i)}|z)$. If $I(x;y|z)>I_\alpha(x;y|z)$, we conclude that $x$ and $y$ as dependent given $z$; otherwise independent. This threshold (tolerance), of mutual Independence is adopted in the forward selection and backward elimination stages as the termination condition. We can interpret the tolerance as the minimum effective quantity of information. In this sense, in the forward ER we are selecting the functions as well as they add a significant quantity of information to the model, while in the backward ER, we discarding functions as well as the information added by them is below the minimum effective quantity, or negligible.

Let $tol$ bel the estimated tolerance of mutual Independence, then we set the termination condition in the backward elimination to:
\begin{equation}
    HLT2 = \begin{cases} 1, \text{ IF } I(\dot{X}; \mathcal{V}(\dot{X},\Phi_i)|\mathcal{V}(\dot{X},\Phi_{\{s_{k-1}-i\}})) > tol, \\ 0, \text{otherwise.} \end{cases}
\end{equation}
where $I(\dot{X}; \mathcal{V}(\dot{X},\Phi_i)|\mathcal{V}(\dot{X},\Phi_{\{s_{k-1}-i\}}))$ follows from Eq.~\ref{eq:backward}.

While in the forward selection stage, we further require the following. Let $I_{a} = I(\dot{X}; \mathcal{V}(\dot{X},\Phi))$ be the mutual information between $\dot{X}$, and the recovered signal $\mathcal{V}(\dot{X},\Phi)$ using \textbf{all} the candidate functions, and let $I_{s} = I(\dot{X}; \mathcal{V}(\dot{X},\Phi_s))$ be the mutual information between $\dot{X}$, and the recovered signal $\mathcal{V}(\dot{X},\Phi_s)$ using the candidate functions indicated by the subset of indices $s$ at the end of each iteration of the forward selection. Then, at any iteration, if $(I_a - I_s) < tol$, that indicates that no other function than what has already been included in the subset $s$ contains significant information if included. So, we formulate the termination condition for the forward selection as:
\begin{equation}
    HLT1 = \begin{cases} 1, \text{ IF } I(\dot{X}; \mathcal{V}(\dot{X},\Phi_i)|\mathcal{V}(\dot{X},\Phi_{s_{k-1}})) < tol, \text{ OR } (I_a - I_s) < tol \\ 0, \text{otherwise.} \end{cases}
\end{equation}
where $I(\dot{X}; \mathcal{V}(\dot{X},\Phi_i)|\mathcal{V}(\dot{X},\Phi_{s_{k-1}}))$ is what indicated in Eq.~\ref{eq:forward}.

\section{Implementation and architecture}
\begin{figure}
    \centering
    \includegraphics[scale=0.55]{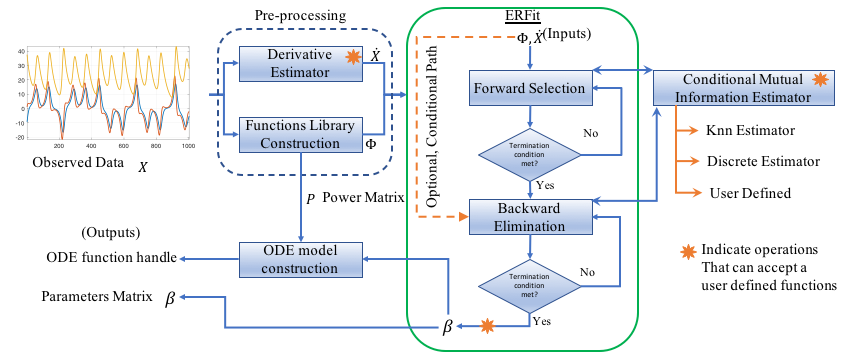}
    \caption{Process flow of the Entropic Regression package (ERFit)}
    \label{fig:erfit}
\end{figure}

Fig.~\ref{fig:erfit} shows the general structure of the ERFit software. We provide a wide range of options and adaptability of our implementation, helping the users (researchers) easily use the package to discover underlying dynamics in science and engineering. Moreover, we make it possible for the researchers to explore and investigate different alternatives by allowing user-defined functions for the mutual information estimator, derivative estimator, and parameter magnitude estimator based on the recovered sparse structure.

Detailed comments and user guide (using MATLAB live script) are provided within the code, to explain every step of the algorithm, and all possible choices for the options.

\section{Quality control and Benchmarking}
The Entropic Regression implementation have been tested on a wide range of systems with varying number of dimensions such as Logistic map (1D), Lorenz (3D), Rossler (3D), coupled network of Logistic map (100D), Kuramoto-Sivenshsky Equations (16D), coupled network of Lorenz system (300D). Note that the dimension mentioned here is the system dimension, which is highly increased after the construction of the candidate functions Library. For example \cite{almomani2019prediction,almomani2020entropic}, the 16D of Kuramoto-Sivashinsky equation, will result with about 150 dimensions for the regression problem, and the 300D of the Lorenz system network (100 nodes, 3D each), considering the second order power polynomial expansion, will result with 45000+ dimensions for the regression problem.

\begin{figure}
    \centering
    \includegraphics[scale=0.6]{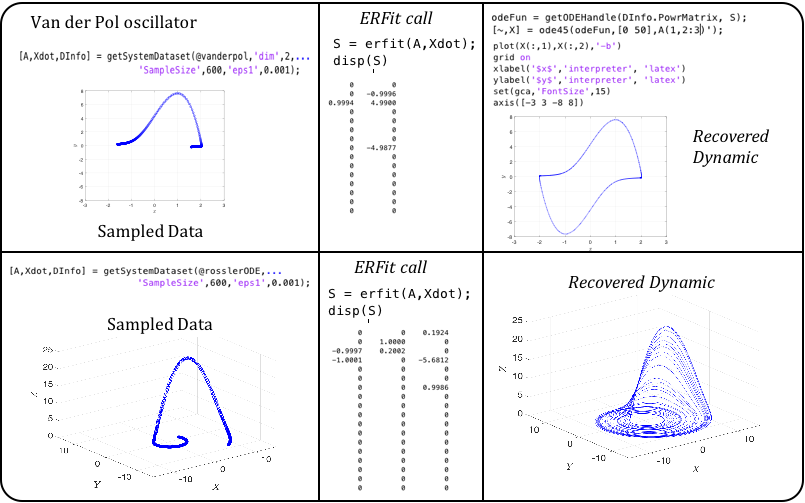}
    \caption{Two illustrative examples for the use and performance of the ERFit.}
    \label{fig:examples}
\end{figure}

In our implementation, we provide a detailed user guide using the ERFit package, and we showed the results for different standard benchmarking problems such as time series from Lorenz system and Rossler attractor. In Fig.~\ref{fig:examples} we present two summarized examples from Rossler and Van der Pol attractors.

\section{Availability}
\vspace{0.5cm}
\section*{Operating system}
The operating system requirements are subject to MATLAB version requirements, which discussed in the Programming language section.

\section*{Programming language}
This code created with MATLAB 2018b. However, for the ERFit package, we emphasize the use of the basic operations while minimal using built-in functions, and according to Matlab release notes (can be found by running the syntax: \texttt{web(fullfile(docroot, 'matlab/release-notes.html'))}, in the command line), the ERFit package is compatible with MATLAB 2016a and newer versions.


\section*{Additional system requirements}
There are no additional requirements. However, for high dimensional problems, it is recommended to have at least 8GB RAM.






\section*{Software location:}



{\bf Code repository:} GitHub: \url{https://github.com/almomaa/ERFit-Package}





\section{Reuse potential}
In a wide range of scientific fields, such as but not limited to biology, epidemiology, chemistry, physics, control systems, and causality inference, a core objective is the data-driven discovery of the underlying dynamics by the construction of a mathematical model that helps to analyze and predict the observed system. The \texttt{erfit} can be used to construct such a reduced dimension mathematical model. It takes the time series observed in a specific experimental setting and produces a reduced dimension mathematical model that is ready to plug-in to the ode solver and build analysis and predictions.  Although the \texttt{erfit} requires minimal supervision and experience from the user, we advance and facilitate the possibility that the researchers explore the different variations with a wide range of options. We made it very simple and straight forward to use a user-defined function for the derivative estimation, mutual information estimator, conditional mutual information estimator, and the parameters magnitude estimator based on the recovered sparse structure.



\section*{Funding statement}
This work was funded in part by the Simons Foundation (Grant No. 318812), the Army Research Office (Grant No. W911NF-16-1-0081), the Office of Naval Research (ONR) (Grant No. N00014-15-1-2093), and the DARPA


\section*{Competing interests}
The authors declare that they have no competing interests.

\bibliography{erfitRef}










\end{document}